# A REMARK ON THE CONVERGING-INPUT CONVERGING-STATE PROPERTY

Eduardo D. Sontag[*]
Department of Mathematics
Rutgers University, New Brunswick, NJ 08903
http://www.math.rutgers.edu/~sontag

**Abstract**

Suppose that an equilibrium is asymptotically stable when external inputs vanish. Then, every bounded trajectory which corresponds to a control which approaches zero and which lies in the domain of attraction of the unforced system, must also converge to the equilibrium. This "well-known" but hard-to-cite fact is proved and slightly generalized here.

## 1 Introduction

This note deals with finite-dimensional controlled systems of the general form

$$\dot{x}(t) = f(x(t), u(t)) \tag{1}$$

and stability properties of an equilibrium. Suppose that $f(0,0) = 0$ and that the zero state is a globally asymptotically stable equilibrium for the "unforced" system $\dot{x} = f(x,0)$. It is well-known that even small inputs $u(\cdot)$ may destabilize the system; in fact there are examples where $u(t) \to 0$ as $t \to \infty$ but $x(t)$ does not converge to zero (or even diverges to infinity).

On the other hand, if $u(t) \to 0$ then *boundedness* of the trajectory $x(\cdot)$ is enough to guarantee $x(t) \to 0$ as $t \to \infty$. This is a "well-known" fact, and one proof was given in [1]. Unfortunately, the conference paper [1] is not easily accessible. In addition, certain (not at all essential) assumptions were made, for simplicity of exposition, which render the result not immediately applicable in some contexts, such as those involving positivity constraints on inputs and states (as in biological and chemical applications). In this note, we basically repeat the proof from that reference, but adapt it to a more general situation, relaxing the global asymptotic stability assumption and allowing inputs to belong to more general spaces than those in [1].

### 1.1 Systems, Notations, Statement of Result

We consider arbitrary finite-dimensional systems of the form (1), where states $x(t)$ take values on an open subset $\mathbb{X}$ of a Euclidean space $\mathbb{R}^n$, for some integer $n$, and inputs $u(t)$ take values on a metric space $\mathbb{U}$. The function

$$f : \mathbb{X} \times \mathbb{U} \to \mathbb{R}^n$$

---

[*]Supported in part by US Air Force Grant F49620-01-1-0063

is continuous, and it satisfies the following local Lipschitz condition: for each $\xi \in \mathbb{X}$ and each compact subset $\mathbb{U}_0$ of $\mathbb{U}$, there exist a neighborhood $\mathcal{V}$ of $\xi$ and a constant $L$ such that $|f(x,u) - f(z,u)| \leq L |x-z|$ (Euclidean norms) for all $x, z \in \mathcal{V}$ and all $u \in \mathbb{U}_0$. (This local Lipschitz condition is satisfied, in particular, if $f(x,u)$ is continuously differentiable on $x$, for each fixed $u$, and its Jacobian $\partial f/\partial x$ is continuous on $\mathbb{X} \times \mathbb{U}$. Below, we will refer to Theorem 1 in [2], which had been stated, purely for reasons of simplicity of exposition, under the assumption that this stronger differentiability condition holds. However, the proof of the result to be quoted relies only upon a theorem given in an appendix to [2] which involves merely the local Lipschitz condition.)

We will assume that two special elements $\bar{x} \in \mathbb{X}$ and $\bar{u} \in \mathbb{U}$ have been singled-out, so that $\bar{x}$ is a steady-state when the input is constantly equal to $\bar{u}$; that is, $f(\bar{x}, \bar{u}) = 0$.

An *input defined on* $\mathcal{I}$, where $\mathcal{I}$ is a subinterval of $[0, \infty)$, is a Lebesgue-measurable function $u : \mathcal{I} \to \mathbb{U}$ which is locally essentially bounded, in the sense that for each compact subset $\mathcal{I}_0$ of $\mathcal{I}$ there is some compact subset $\mathbb{U}_0 \subseteq \mathbb{U}$ such that $u(t) \in \mathbb{U}_0$ for almost all $t \in \mathcal{I}_0$. Given any input $u$ defined on an interval $\mathcal{I}$ containing 0, and any initial state $\xi$, there is a unique maximal solution $x(t) = \varphi(t, \xi, u)$ of (1) with initial value $x(0) = \xi$; this solution is defined on some maximal interval $[0, \sigma_{\xi, u})$ of $\mathcal{I}$.

When $u(t) \equiv \bar{u}$, we write $\varphi(t, \xi, u)$ simply as $\varphi(t, \xi)$; this is the solution $x(t)$ of the autonomous system $\dot{x} = f(x, \bar{u})$ with $x(0) = \xi$. Note that $\varphi(t, \bar{x}) = \bar{x}$ for all $t \geq 0$.

From now on, we will denote by single bars "$|\cdot|$" the distances to the "origins" $\bar{x}$ in $\mathbb{X}$ or $\bar{u}$ in $\mathbb{U}$: $|\xi| = \text{dist}(\xi, \bar{x})$, $|\mu| = \text{dist}(\mu, \bar{u})$ and we will use double bars "$\|\cdot\|$" for the supremum norm on the spaces of controls and of trajectories. That is, if $x : \mathcal{I} \to \mathbb{X}$ is an absolutely continuous function and $u : \mathcal{I} \to \mathbb{X}$ is an input, where $\mathcal{I} \subseteq [0, \infty)$ is an interval, then

$$\|x\| = \sup_{t \in \mathcal{I}} |x(t)| \quad \text{and} \quad \|u\| = \sup_{t \in \mathcal{I}} |u(t)|$$

where the second "sup" is interpreted as an essential supremum.

One last item of terminology is as follows. Given a compact subset $K \subseteq \mathbb{X}$, we will say that a function $x : \mathcal{I} \to \mathbb{X}$ is *K-recurrent* if for each $T > 0$ there is some $t > T$ such that $x(t) \in K$. (A weaker notion would result if asking merely that $\Omega^+[x] \bigcap K \neq \emptyset$, where $\Omega^+[x]$ is the omega-limit set of $x$. However, this property amounts simply to $K'$-recurrence for a slightly larger compact set $K \subseteq K'$.)

This is the result that we wish to prove:

**Theorem 1** *Suppose that $\bar{x}$ is an asymptotically stable equilibrium of the autonomous system $\dot{x} = f(x, \bar{u})$, with domain of attraction $\mathcal{O}$, and that $K$ is a compact subset of $\mathcal{O}$. Let $x(\cdot)$ be a K-recurrent solution of (1) defined on $[0, \infty)$, and suppose that $u(t) \to \bar{u}$ as $t \to \infty$. Then, $x(t) \to \bar{x}$ as $t \to \infty$.*

*Furthermore, the following stability property holds: for each $\varepsilon > 0$ there is some $\delta > 0$ such that, whenever $|\xi| < \delta$ and $\|u\| < \delta$, the solution $x(t) = \varphi(t, \xi, u)$ exists for all $t \geq 0$, and $|x(t)| < \varepsilon$ for all $t \geq 0$.*

Recall that the domain of attraction $\mathcal{O}$ of an asymptotically stable equilibrium $\bar{x}$ of the autonomous system $\dot{x} = f(x, \bar{u})$ is the set consisting of those initial conditions $\xi$ for which $\varphi(t, \xi) \to \bar{x}$ as $t \to \infty$; the set $\mathcal{O}$ is an open subset of $\mathbb{X}$.

We remark that a system which satisfies the above properties: stability and $x(t) \to 0$ whenever $u(t) \to 0$, is not necessarly an ISS system. A counterexample is $\dot{x} = (-1+u)x$ with $\mathbb{X} = \mathbb{R}$ and $\mathbb{U} = [0,1]$.

## 2 Proofs

We collect first a number of statements, all of which are elementary consequences of continuity properties of solutions of differential equations. We assume given a system (1) so that $\bar{x}$ is an asymptotically stable equilibrium of the autonomous system $\dot{x} = f(x, \bar{u})$, with domain of attraction $\mathcal{O}$.

Asymptotic stability implies in particular stability, i.e., for each $\varepsilon > 0$ there is some $\delta = \Delta_1(\varepsilon) > 0$ such that:

$$|\xi| \leq \delta \quad \Rightarrow \quad |\varphi(t,\xi)| < \varepsilon \quad \forall \, t \geq 0 \,. \tag{2}$$

**Lemma 2.1** For each compact subset $K \subseteq \mathcal{O}$, each $T \geq 0$, and each $\varepsilon > 0$, there is a $\delta = \Delta_2(K, T, \varepsilon) > 0$ with the following property: for each $\xi \in K$ and every control $u$ defined on $[0, \infty)$ such that $\|u\| < \delta$, the solution $\varphi(t, \xi, u)$ is defined for $t \in [0, T]$, and

$$\operatorname{dist}\left(\varphi(t, \xi, u), \varphi(t, \xi)\right) < \varepsilon \quad \forall \, t \in [0, T] \,. \tag{3}$$

*Proof.* We recall Theorem 1 from [2]. For each fixed $T > 0$, let $\mathcal{D}_T$ be the set consisting of those pairs $(\xi, u)$ with $u : [0, T] \to \mathbb{U}$, for which $\varphi(t, \xi, u)$ is defined for all $t \in [0, T]$. Then $\mathcal{D}_T$ is an open subset of $\mathbb{X} \times \mathbb{U}_T$, where $\mathbb{U}_T = \mathcal{L}_\mathbb{U}^\infty(0, T)$, the space of all Lebesgue measurable essentially bounded $u : [0, T] \to \mathbb{U}$ endowed with the (essential) supremum norm. Moreover, the mapping $(\xi, u) \mapsto x = \varphi(\cdot, \xi, u)$ is continuous, when trajectories are also endowed with the supremum (uniform convergence) norm. Now suppose that $K, T, \varepsilon$ are given. Pick any $\zeta \in K$. By continuity of $\alpha$ at the pair $(\zeta, 0)$, there is some $\delta_\zeta > 0$ such that $\operatorname{dist}(\varphi(t, \xi, u), \varphi(t, \xi', u')) < \varepsilon$ for all $t \in [0, T]$ holds for all states $\xi, \xi'$ and inputs $u, u'$ in balls of radius $\delta_\zeta$ around $\zeta$ and the zero input respectively. In particular, $\operatorname{dist}(\varphi(t, \xi, u), \varphi(t, \xi)) < \varepsilon$ for all $t \in [0, T]$. By compactness of $K$, one may let $\delta$ be obtained as the smallest $\delta_\zeta$ from a finite subcover. ∎

Lemma 5.9.12 in [2] states that, for the autonomous system $\dot{x} = f(x, \bar{u})$, for each $\xi \in \mathcal{O}$ and each neighborhood $\mathcal{V}$ of $\bar{x}$ there is some neighborhood $\mathcal{W}$ of $\xi$ and some $T \geq 0$ such that

$$\zeta \in \mathcal{W} \ \& \ t \geq T \quad \Rightarrow \quad \varphi(t, \zeta) \in \mathcal{V} \,.$$

A compactness argument then gives the following standard uniform stability fact.

**Lemma 2.2** For each compact subset $K \subseteq \mathcal{O}$, and each neighborhood $\mathcal{V}$ of $\bar{x}$, there is some $T = \mathcal{T}(K, \mathcal{V}) > 0$ such that $\varphi(T, \xi) \in \mathcal{V}$ for all $\xi \in K$. □

**Lemma 2.3** For each compact subset $K \subseteq \mathcal{O}$ and each $\varepsilon > 0$, there exist $T = \mathcal{T}(K, \varepsilon)$ and $\delta = \Delta_3(K, \varepsilon)$ such that, for every $\xi \in K$ and every input $u$ with $\|u\| < \delta$, the following properties hold:

(a) The solution $\varphi(t, \xi, u)$ is defined for all $t \in [0, T]$.

(b) Denoting $\widetilde{K} = \{\varphi(t,\xi) \mid t \in [0,T], \xi \in K\}$, it holds that $\text{dist}\,(\varphi(t,\xi,u),\widetilde{K}) < \varepsilon$ for all $t \in [0,T]$.

  (c) $|x(T,\xi,u)| < \varepsilon$.

*Proof.* Given $K,\varepsilon$, introduce $\mathcal{V} =$ the open ball of radius $\varepsilon/2$ around $\bar{x}$, $T = \mathcal{T}(K,\varepsilon)$ chosen as in Lemma 2.2, and $\delta = \Delta_2(K,T,\varepsilon/2)$ where $\Delta_2$ is the function in Lemma 2.1. Now pick any $\xi \in K$ and any $\|u\| < \delta$. Inequality (3) says that $\text{dist}\,(\varphi(t,\xi,u),\varphi(t,\xi)) < \varepsilon/2$ for all $t \in [0,T]$, which gives in particular property (b). On the other hand, this inequality together with $|\varphi(T,\xi)| < \varepsilon/2$ (which holds because of the choice of $T$) gives property (c). ∎

Next, we prove the stability part of Theorem 1.

Let $\varepsilon > 0$ be given. We introduce $\delta_1 = \min\{\Delta_1(\varepsilon/2),\varepsilon\}$, $K =$ the closed ball around $\bar{x}$ of radius $\delta_1/2$, and $T = \mathcal{T}(K,\delta_1/2)$ and $\delta = \Delta_3(K,\delta_1/2)$ as in Lemma 2.3. We prove the stability property with this $\delta$.

Let $x \in K$ and $\|u\| < \delta$. By Lemma 2.3, property (c), $|x(T,\xi,u)| < \delta_1/2$, so
$$\xi' = x(T,\xi,u) \in K\,.$$

On the other hand, the choice $\delta_1 \leq \Delta_1(\varepsilon/2)$ insures $|\varphi(t,\xi)| < \varepsilon/2$ for all $t$, and in particular for $t \in [0,T]$, so $\widetilde{K}$ is included in the $\varepsilon/2$-ball around $\bar{x}$, so the triangle inequality applied to
$$\text{dist}\,(\varphi(t,\xi,u),\widetilde{K}) < \frac{\delta_1}{2} < \frac{\varepsilon}{2} \quad \forall\, t \in [0,T]$$
(property (b) in Lemma 2.3) lets us conclude that
$$|\varphi(t,\xi,u)| < \varepsilon \quad \forall\, t \in [0,T]\,.$$

Starting now from the state $\xi'$, and using the restriction of the input $u$ to $[T,\infty)$, which also has uniform norm $< \delta$, we conclude that $|\varphi(t,\xi,u)| < \varepsilon$ for all $t \in [T,2T]$, and an induction argument proves that this holds for all $t \geq 0$ as required for the stability proof.

To prove the first part of the Theorem, we first show the following result.

**Proposition 2.4** For each compact $K \subseteq \mathcal{O}$ and each $\varepsilon > 0$ there exist $T_0 \geq 0$ and $\delta > 0$ with the following properties: for every control $u$ defined on $[0,\infty)$ such that $\|u\| < \delta$ and every $\xi \in K$, the solution $x(t) = \varphi(t,\xi,u)$ is defined for all $t \geq 0$ and satisfies that $|x(t)| \leq \varepsilon$ for all $t \geq T_0$.

*Proof.* Pick any compact $K$ and $\varepsilon > 0$. Let $\delta_1 > 0$ be such that the stability statement holds, i.e. $|\varphi(t,\xi,u)| < \varepsilon$ for all $t \geq 0$ provided that $|\xi| < \delta_1$ and $\|u\| < \delta_1$. Pick $T_0 = \mathcal{T}(K,\delta_1)$ and $\delta_2 = \Delta_3(K,\delta_1)$ as in Lemma 2.3, and let $\delta = \min\{\delta_1,\delta_2\}$. Now take any $u$ with $\|u\| < \delta$ and any $x \in K$. By Lemma 2.3, and using that $\|u\| < \delta_2$, the solution is defined on $[0,T_0]$ and $|\xi'| < \delta_1$, where $\xi' = \varphi(T_0,\xi,u)$. Since the restriction $u'$ of $u$ to $[T_0,\infty)$ has norm $< \delta_1$, the stability statement applied to the initial state $\xi'$ and input $u'$ insures that $|\varphi(t,\xi,u)| < \varepsilon$ for all $t \geq T_0$. ∎

Finally, we prove the first part of the Theorem. Let $K$ be so that $x(t) = \varphi(t,\xi,u)$ is $K$-recurrent and suppose that $u(t) \to \bar{u}$ as $t \to \infty$. We must show that for each $\varepsilon > 0$ there is some $T > 0$ such that $|x(t)| < \varepsilon$ for all $t \geq T$. Pick $T_0$ and $\delta$ as in Proposition 2.4. Since $u(t) \to \bar{u}$, we may pick some $T_1$ such that $\|u\| < \delta$ for all $t \geq T_1$. By the $K$-recurrence property, we may pick $T_2 \geq T_1$ such that $\xi' = \varphi(T_2, \xi, u)) \in K$. Starting from the initial state $\xi'$ and using the input $u$ restricted to $[T_2, \infty)$, we are in the situation of Proposition 2.4, and this insures that $|\varphi(t, \xi, u)| < \varepsilon$ for all $t \geq T := T_2 + T_0$. ∎

## Acknowledgment

The author wishes to thank David Angeli for the suggestion that the main result should be stated in terms of a recurrence property.

## References

[1] Sontag, E.D., "Remarks on stabilization and input-to-state stability," *Proc. IEEE Conf. Decision and Control, Tampa, Dec. 1989*, IEEE Publications, 1989, pp. 1376-1378.

[2] Sontag, E.D., *Mathematical Control Theory: Deterministic Finite Dimensional Systems*, Springer, New York, 1990. Second Edition, 1998.